\documentclass[12pt,twoside]{article}
\usepackage{amsfonts}
\usepackage[dvips]{graphicx}
\usepackage{amsmath,fleqn}
\usepackage{amssymb}
\usepackage{mathrsfs}
\usepackage{amstext}

\numberwithin{figure}{section} \numberwithin{equation}{section}
\makeatletter 
\author{\small{ Lin Zhao\footnote{Corresponding author. Fax: +86-516-83591530. E-mail address: zhaolinmath@163.com (L. Zhao) Research supported by the Natural Science Foundation of Jiangsu Province of China (No. 20180638)} }
\medskip\\
\emph{\small{Department of Mathematics,
China University of Mining and Technology}}\\
\emph{\small{ Xuzhou, Jiangsu 221116, China}  }}
\begin{document}
\date{}
\title{\textbf{Cocompact imbedding theorem for functions of bounded variation into Lorentz spaces }} \maketitle

\noindent{\bf Abstract}
\medskip\\
\indent \small{We show that the imbedding $\dot{BV}(\mathbb{R}^N)\hookrightarrow L^{1^\ast,q}(\mathbb{R}^N)$, $q>1$ is cocompact with respect to group and the profile decomposition for $\dot{BV}(\mathbb{R}^N)$. This paper extends the cocompactness and profile decomposition for the critical space $L^{1^\ast}(\mathbb{R}^N)$ to Lorentz spaces $L^{1^\ast,q}(\mathbb{R}^N)$, $q>1$. A counterexample for $\dot{BV}(\mathbb{R}^N)\hookrightarrow L^{1^\ast,1}(\mathbb{R}^N)$ not cocompact is given in the last section. }
\medskip\\
\section{Introduction and Preliminaries}
It is well know that the classical Sobolev imbedding may be improved within the framework of Lorentz spaces $L^{p,q}$, see \cite{CRT,Pe,SS,LT1}. Solimini \cite{SS} shows that the space $W^{1,p}(\mathbb{R}^N)$, $1<p<N$, imbeds into $L^{p^\ast,q}(\mathbb{R}^N)$, $p\leq q\leq +\infty$ and give the profile decomposition for Sobolev space $W^{1,p}(\mathbb{R}^N)$. Adimurthi and Tintarev \cite{AC} extends the profile decomposition proved by Solimini \cite{SS} for Sobolev spaces $W^{1,p}(\mathbb{R}^N)$ with $1<p<N$ to the non-reflexive case $p=1$. They replace $W^{1,1}(\mathbb{R}^N)$ by $\dot{BV}(\mathbb{R}^N)$ and show the imbedding $\dot{BV}(\mathbb{R}^N)\hookrightarrow L^{1^\ast}(\mathbb{R}^N)$ is cocompact. In this paper
we extend the result of \cite{AC} to the framework of Lorentz spaces and give the cocompact imbedding and the profile decomposition for Lorentz spaces. In broad sense defect of compactness is known as concentration compactness by Lions \cite{L1}. It is expressed as a sum of elementary concentrations called profile decomposition, which is introduced by Struwe \cite{MS}, also see Lieb \cite{HL}. G\'{e}rard \cite{Ge} and Jaffard \cite{Ja} study the profile decomposition for the fractional Sobolev spaces and give a weaker form of remainder. Schihdler and Tintarev \cite{ST2} prove the profile decomposition for general Hilbert spaces. Profile decomposition for uniformly convex Banach space is obtained by Solimini and Tintarev \cite{ST,ST3}. A profile decomposition in the general non-reflexive Banach space remain an open problem, see \cite{AC,ST}. The main results of this paper are Theorem 2.1 and 2.4. The idea for the proof of the theorem is followed by \cite{AC,TF}, but uses different argument for the estimate of Lorentz seminorms over lattices.

We need the following definitions and lemmas of Lorentz spaces and $BV$ spaces, see \cite{ATL,G,NI,LT1,LT2,PZ}. Lorentz spaces are known as real interpolation spaces between Lebesgue spaces and can be defined via the notion of Schwarz symmetrization. Let $u$ be a measurable function on $\mathbb{R}^N$, $N\geq2$, whose level sets have finite measure for every level. Then the function
\begin{eqnarray}\label{}
\mu(\lambda)=|\{x\in \mathbb{R}^N\mid|u(x)|>\lambda\}|,~\lambda\geq0,
\end{eqnarray}
is the distribution function of $u$ and
\begin{eqnarray}\label{}
u^\sharp(r)=\inf\{\lambda>0\mid \mu(\lambda)\leq |B_r|\}
\end{eqnarray}
is radially symmetric and non-increasing rearrangement of $u$, where $B_r$ is ball centered at the origin with radius $r$ and $|A|$ is the $N$ dimensional Lebesgue measure of $A\subset\mathbb{R}^N$. We call the function $u^\sharp(|x|)$, $x\in \mathbb{R}^N$ as Schwarz symmetrization of $u$. Then, we define Lorentz spaces $L^{p,q}(\mathbb{R}^N)$ as
\begin{eqnarray}\label{1.3}
L^{p,q}(\mathbb{R}^N)=\{u:\text{measurable in}~\mathbb{R}^N\mid \|u\|_{L^{p,q}}=|B_1|^{\frac{q-p}{pq}}(\int_{\mathbb{R}^N}(|x|^\frac{N}{p}u^\sharp(|x|))^q\frac{dx}{|x|^N})^\frac{1}{q}<\infty\}.\nonumber\\
\end{eqnarray}
The definition of Lorentz spaces is also given by decreasing rearrangement $u^\ast$ as the distribution function $\mu$ of $u$, namely
\begin{eqnarray}\label{}
u^\ast(t)=\sup\{\lambda>0\mid \mu(\lambda)>t\}.
\end{eqnarray}
Then, we define the Lorentz spaces
\begin{eqnarray}\label{1.5}
L^{p,q}(\mathbb{R}^N)=\{u:\text{measurable in}~\mathbb{R}^N\mid\|u\|_{L^{p,q}}=(\int_0^\infty(u^\ast(t)t^\frac{1}{p})^q\frac{dt}{t}    )^\frac{1}{q}<+\infty\},
\end{eqnarray}
where (\ref{1.5}) is a quasinorm and equivalent to (\ref{1.3}). Lorentz spaces satisfy the following imbedding lemma.\\
\textbf{Lemma 1.1} (\cite{G}). Suppose that $(\Omega,\Sigma,m)$ is non-atomic and that $m(\Omega)<\infty$. Then if $0<p_1<p_2\leq \infty$, $L^{p_2,q_2}(\Omega)\subset L^{p_1,q_1}(\Omega)$ for any $q_1$ and $q_2$, with continuous inclusion, where $m(A)$ is the $N$ dimensional Lebesgue measure of $A\subset\mathbb{R}^N$.\\
\textbf{Definition 1.2} (\cite{PZ}). Let $\Omega$ be an open subset of $\mathbb{R}^N$. A function $u\in L^1(\Omega)$ is said to be of bounded variations, if the total variation of $u$ on $\Omega$ is
\begin{eqnarray}\label{}
\|Du\|_{\Omega}=\sup \{\int_\Omega u\text{div} v dx \mid v \in C_0^\infty (\Omega;\mathbb{R}^N),~\|\varphi\|_\infty \leq 1\}<+\infty,
\end{eqnarray}
where $v(x)=(v_1(x),v_2(x),...v_N(x))$ and $\|v\|_\infty =\sup_{x\in \Omega}(\sum_{k=1}^{N}(v_k(x))^2)^\frac{1}{2}$. On
\begin{eqnarray}\label{}
BV(\Omega)=\{ u\in L^1(\Omega) \mid \|Du\|_{\Omega}<\infty\},
\end{eqnarray}
we define the norm
\begin{eqnarray}\label{}
\|u\|_{BV(\Omega)}=\|Du\|_\Omega +\|u\|_{L^1(\Omega)}.
\end{eqnarray}
\textbf{Definition 1.3.} The sequence $\{u_n\}$ converges weakly to $u$ in $BV(\Omega)$, written $u_n\rightharpoonup u$ if
\begin{eqnarray}\label{}
u_n\rightarrow u~\text{in}~ L^1_{loc}(\Omega),
\end{eqnarray}
and
\begin{eqnarray}\label{}
\partial_k u_n\rightharpoonup\partial_k u~\text{in}~(C_0(\Omega))^\ast,
\end{eqnarray}
for $1\leq k\leq N$, as $n\rightarrow\infty$, where $(C_0(\Omega))^\ast$ denotes the space of finite measures on $\Omega$.\\
\textbf{Lemma 1.4.} Assume $\Omega\subset \mathbb{R}^N$ is open and bounded with $\partial\Omega$ Lipschitz. The imbedding
\begin{eqnarray}\label{}
BV(\Omega)\hookrightarrow L^{1^\ast,1}(\Omega)\subset L^{1^\ast,q}(\Omega),~1<q\leq+\infty,
\end{eqnarray}
is continuous. More precisely, there exists a constant $C$ with depends only on $\Omega,q$ and $N$, such that for $\forall$ $u\in BV(\Omega)$
\begin{eqnarray}\label{}
\|u\|_{L^{1^\ast,q}(\Omega)}\leq\|u\|_{L^{1^\ast,1}(\Omega)}\leq C\|u\|_{BV(\Omega)}.
\end{eqnarray}
\textbf{Definition 1.5.} The space of functions of bounded variation $\dot{BV}(\mathbb{R}^N)$ is the space of all measurable functions $u$ vanishing at infinity (i.e. $\forall$ $M>0$, $|\{x\in \mathbb{R}^N\mid|u(x)|>M\}|<\infty$) such that
\begin{eqnarray}\label{}
\|u\|_{\dot{BV}}=\|Du\|=\sup\{\int_{\mathbb{R}^N} u\text{div} v dx \mid v\in C_0^\infty(\mathbb{R}^N;\mathbb{R}^N),\|v\|_{L^\infty}\leq1 \}<\infty,
\end{eqnarray}
where $v(x)=(v_1(x),v_2(x),...v_N(x))$ and $\|v\|_\infty =\sup_{x\in \Omega}(\sum_{k=1}^{N}(v_k(x))^2)^\frac{1}{2}$. The $\dot{BV}(\mathbb{R}^N)$ norm $\|u\|_{\dot{BV}}$ can be interpreted as the total variation $\|Du\|$ of the measure associated with to derivative $Du$ in the sense of distributions on $\mathbb{R}^N$.\\
\textbf{Definition 1.6.} The sequence $\{u_n\}$ converges weakly to $u$ in $\dot{BV}(\mathbb{R}^N)$, written $u_n\rightharpoonup u$ if
\begin{eqnarray}\label{}
u_n\rightarrow u~\text{in}~L^1_{loc}(\mathbb{R}^N),
\end{eqnarray}
and
\begin{eqnarray}\label{}
\partial_k u_n\rightharpoonup\partial_k u,~1\leq k\leq N
\end{eqnarray}
as finite measures on $\mathbb{R}^N$ and $n\rightarrow\infty$.\\
\textbf{Lemma 1.7} (\cite{LT2}). $\dot{BV}(\mathbb{R}^N)$ is continuously embedding in $L^{1^\ast,1}(\mathbb{R}^N)$ for $N\geq2$, and in $L^\infty(\mathbb{R})$ for $N=1$.

We need the following properties of $\dot{BV}(\mathbb{R}^N)$, also see \cite{AC}.

1). Invariance. The group of operators on $\dot{BV}(\mathbb{R}^N)$,
\begin{eqnarray}\label{1.16}
\mathcal{D}=\{g[j,y]:u\rightarrow2^{(N-1)j}u(2^j(\cdot-y))         \}_{j\in \mathbb{Z}},
\end{eqnarray}
consists for linear isometries of $\dot{BV}(\mathbb{R}^N)$, which are also linear isometries on $L^{1^\ast,q}(\mathbb{R}^N)$, $1\leq q\leq \infty$.

2). Chain rule. Let $\varphi\in C^1(\mathbb{R})$. Then for $\forall$ $u\in \dot{BV}(\mathbb{R}^N)$
\begin{eqnarray}\label{1.17}
\|D\varphi(u)\|\leq\|\varphi\|_{L^{\infty}} \|Du\|.
\end{eqnarray}
\textbf{Definition 1.8} (\cite{AC,TF}). Let $X$  be a Banach space and let $\mathcal{D}$ be a group of linear isometries of $X$. One says that a sequence $\{u_k\}\subset X$ is $\mathcal{D}$-vanishing (to be written $u_k\stackrel{\mathcal{D}}{\rightharpoonup}0$) if for any sequence $\{g_k\}\subset \mathcal{D}$ one has $g_ku_k\rightharpoonup 0$ in $X$. A continuous imbedding of $X$ into a topological space $Y$ is called cocompact with respect to $\mathcal{D}$ if $u_k\stackrel{\mathcal{D}}{\rightharpoonup}0$ implies $u_k\rightarrow0$ in $Y$.

\section{Main results}
\textbf{Theorem 2.1.} The imbedding $\dot{BV}(\mathbb{R}^N)\hookrightarrow L^{1^\ast, q}(\mathbb{R}^N)$, $1<q\leq 1^\ast=\frac{N}{N-1}$, $N\geq2$, is cocompact with respect to the group $\mathcal{D}$ (\ref{1.16}), i.e. if, for any sequence $\{j_k,y_k\}\subset \mathbb{Z}\times \mathbb{R}^N$, $g[j_k,y_k]u_k\rightharpoonup 0$ in $\dot{BV}(\mathbb{R}^N)$ then $u_k\rightarrow0$ in $L^{1^\ast, q}(\mathbb{R}^N)$.\\
\textbf{Remark 2.2.} The embedding $\dot{BV}(\mathbb{R}^N)\hookrightarrow  L^{1^\ast,q}(\mathbb{R}^N)$, $1^\ast<q\leq\infty$ is cocompact with respect to $\mathcal{D}$, by the cocompact imbedding $\dot{BV}(\mathbb{R}^N)\hookrightarrow  L^{1^\ast}(\mathbb{R}^N)$ and $L^{1^\ast}(\mathbb{R}^N)\subset L^{1^\ast,q}(\mathbb{R}^N)$, $q>1^\ast$. Theorem 2.1 may hold as follows by the interpolate method, since the sequence is bounded in $L^{1^\ast,1}(\mathbb{R}^N)$, and we will use the estimate of Lorentz norms over lattices to prove Theorem 2.1. However, the cocompact embedding $\dot{BV}(\mathbb{R}^N)\hookrightarrow L^{1^\ast,1}(\mathbb{R}^N)$ does not hold for $q=1$ and we will give a counterexample in the Section 3 for this case. \\
\textbf{Proof.} Let $\{u_k\}\subset \dot{BV}(\mathbb{R}^N)$ be such that for any $\{j_k,y_k\}\subset \mathbb{Z}\times \mathbb{R}^N$,
\begin{eqnarray}\label{}
g_ku_k=g[j_k,y_k]u_k\rightharpoonup 0~\text{in}~\dot{BV}(\mathbb{R}^N).
\end{eqnarray}
\textbf{Step 1.} Assume first that $\sup_{k\in\mathbb{N}}\|u_k\|_{L^\infty}<\infty$ and $\sup_{k\in\mathbb{N}}\|u_k\|_{L^1}<\infty$. Then using the $L^\infty$-boundedness of $\{u_k\}$ and $BV(\Omega_0)\hookrightarrow L^{1^\ast,q}(\Omega_0)$, $\Omega_0=(0,1)^N$, $1<q\leq 1^\ast=\frac{N}{N-1}$, we have
\begin{eqnarray}\label{}
\|u_k\|_{L^{1^{\ast},q}(\Omega_0)}&=&|B(0,1)|^\frac{q-1^\ast}{1^\ast q}(\int_{\Omega_0}(u_k^\sharp(|x|))^q |x|^{N(q-1)-q}dx           )^\frac{1}{q} \nonumber\\
&\leq& C\|u_k\|_{BV(\Omega_0)}\nonumber\\\
&=&C(\|Du_k\|_{\Omega_0}+\|u_k\|_{L^1(\Omega_0)}),
\end{eqnarray}
which implies
\begin{eqnarray}\label{2.3}
&&|B(0,1)|^\frac{q-1^\ast}{1^\ast q}\int_{\Omega_0}(u_k^\sharp(|x|))^q |x|^{N(q-1)-q}dx \nonumber\\\
 &\leq& C(\|Du_k\|_{\Omega_0}+\|u_k\|_{L^1(\Omega_0)})(\int_{\Omega_0}(u_k^\sharp(|x|))^q |x|^{N(q-1)-q}dx           )^\frac{q-1}{q}.
\end{eqnarray}
Considering Lemma 1.1 and setting $p_1=1^\ast=\frac{N}{N-1}<p_2=3$, $q_1=q$, $q_2=3$ and $m(\Omega_0)=m((0,1)^N)<\infty$, we have
\begin{eqnarray}\label{}
L^3(\Omega_0)=L^{3,3}(\Omega_0)\subset L^{1^\ast,q}(\Omega_0),
\end{eqnarray}
and
\begin{eqnarray}\label{2.5}
\|u_k\|_{L^{1^\ast,q}(\Omega_0)}\leq \|u_k\|_{L^3(\Omega_0)}.
\end{eqnarray}
By (\ref{2.3}) and (\ref{2.5}), we get
\begin{eqnarray}\label{2.6}
&&|B(0,1)|^\frac{q-1^\ast}{1^\ast q}\int_{\Omega_0}(u_k^\sharp(|x|))^q |x|^{N(q-1)-q}dx \nonumber  \\ \nonumber
 &\leq& C(\|Du_k\|_{\Omega_0}+\|u_k\|_{L^1(\Omega_0)})(\int_{\Omega_0}(u_k^\sharp(|x|))^q |x|^{N(q-1)-q}dx           )^\frac{q-1}{q}\\ \nonumber
 &\leq& C(\|Du_k\|_{\Omega_0}+\|u_k\|_{L^1(\Omega_0)})(\int_{\Omega_0} |u_k|^3dx   )^\frac{q-1}{3}\\
 &\leq& C'(\|Du_k\|_{\Omega_0}+\|u_k\|_{L^1(\Omega_0)})(\int_{\Omega_0} |u_k|dx   )^\frac{q-1}{3},
\end{eqnarray}
where the last inequality is given by $L^\infty$-boundedness of $\{u_k\}$. Repeating this inequality (\ref{2.6}) for the domain of intergration $(0,1)^N+y$, $y\in \mathbb{Z}^N$, and adding the resulting inequalities over all $y\in\mathbb{Z}^N$, we have
\begin{eqnarray}\label{}
&&|B(0,1)|^\frac{q-1^\ast}{1^\ast q}\int_{\mathbb{R}^N}(u_k^\sharp(|x|))^q |x|^{N(q-1)-q}dx \nonumber \\
&=&|B(0,1)|^{\frac{q-1^\ast}{1^\ast}(\frac{1}{q}-1)}\|u_k\|^q_{L^{1^\ast,q}(\mathbb{R}^N)} \nonumber \\
&\leq& C''(\|Du_k\|_{\mathbb{R}^N}+\|u_k\|_{L^1(\mathbb{R}^N)})(\sup_{y\in\mathbb{Z}^N}\int_{\Omega_0} |u_k(x-y)|dx   )^\frac{q-1}{3}.
\end{eqnarray}
Here we use the fact that the sum $\sum_{y\in\mathbb{Z}^N}\|Du_k\|_{(0,1)^N+y}$ can be split into $3^N$ sums of variations over unions of cubes with disjoint closures, each of them, as follows from Definition 1.5 bound by $\|Du_k\|_{\mathbb{R}^N}$, which implies $\sum_{y\in\mathbb{Z}^N}\|Du_k\|_{(0,1)^N+y}\leq 3^N\|Du_k\|_{\mathbb{R}^N}$. By the assumption $g[j_k,y_k]u_k\rightharpoonup 0$ in $\dot{BV}(\mathbb{R}^N)$, we have $u_k(\cdot-y_k)\rightarrow 0$ in $L^1((0,1)^N)$ for $\forall$ $\{y_k\}\subset \mathbb{R}^N$. This implies $u_k\rightarrow0$ in $L^{1^\ast,q}(\mathbb{R}^N)$.\\
\\
\textbf{Step 2.} Now consider a general sequence $\{u_k\}\subset \dot{BV}(\mathbb{R}^N)$ satisfying $g[j_k,y_k]u_k\rightharpoonup0$ in $\dot{BV}(\mathbb{R}^N)$ for any $\{j_k,y_k\}\subset \mathbb{Z}\times\mathbb{R}^N$. Let $\chi\in C_0^\infty ((\frac{1}{2^{N-1}}, 4^{N-1}))$ be such that $\chi(t)=t$ if $t\in[1,2^{N-1}]$. Let $\chi_j(t)=2^{(N-1)j}\chi(2^{-(N-1)j}|t|)$, $j\in\mathbb{Z}$, $t\in\mathbb{R}$ and obviously $\|\chi_j'\|_{L^\infty}=\|\chi'\|_{L^\infty}$. By Lemma 1.7 the imbedding $\dot{BV}(\mathbb{R}^N)\hookrightarrow  L^{1^\ast,1}(\mathbb{R}^N) \subset L^{1^\ast,q}(\mathbb{R}^N)$ and $\chi_j(u_k)\in \dot{BV}(\mathbb{R}^N)$, we have
\begin{eqnarray}\label{2.8}
\|\chi_j(u_k)\|_{L^{1^\ast,q}(\mathbb{R}^N)}\leq C\|\chi_j(u_k)\|_{\dot{BV}}=C\|D\chi_j(u_k)\|=C\|D\chi_j(u_k)\|_{B_{kj}}.
\end{eqnarray}
From the definition of $L^{1^\ast,q}(\mathbb{R}^N)$, we get
\begin{eqnarray}\label{2.9}
&&\|\chi_j(u_k)\|_{L^{1^{\ast},q}(\mathbb{R}^N)}=|B(0,1)|^\frac{q-1^\ast}{1^\ast q}(\int_{\mathbb{R}^N}((\chi_j (u_k))^\sharp(|x|))^q |x|^{N(q-1)-q}dx           )^\frac{1}{q} \nonumber  \\ \nonumber
&=& |B(0,1)|^\frac{q-1^\ast}{1^\ast q}(\int_{\{x\in \mathbb{R}^N\mid 0\leq |u_k(x)|<\infty\}       }((\chi_j ( u_k))^\sharp(|x|))^q |x|^{N(q-1)-q}dx           )^\frac{1}{q}\\  \nonumber
&\geq&|B(0,1)|^\frac{q-1^\ast}{1^\ast q}(\int_{\{x \in \mathbb{R}^N \mid 2^{(N-1)j}\leq |u_k(x)|<2^{(N-1)(j+1)}\}       }((\chi_j (u_k))^\sharp(|x|))^q |x|^{N(q-1)-q}dx           )^\frac{1}{q}\nonumber\\
&=& |B(0,1)|^\frac{q-1^\ast}{1^\ast q}(\int_{ A_{kj}       }((\chi_j ( u_k))^\sharp(|x|))^q |x|^{N(q-1)-q}dx           )^\frac{1}{q}
\end{eqnarray}
where $A_{kj}=\{x \in \mathbb{R}^N \mid 2^{(N-1)j}\leq |u_k(x)|<2^{(N-1)(j+1)}\}$, $j\in \mathbb{Z}^N$. Considering the inequalities (\ref{2.8}) and (\ref{2.9}), we have
\begin{eqnarray}\label{}
|B(0,1)|^\frac{q-1^\ast}{1^\ast q}(\int_{ A_{kj}       }((\chi_j(u_k))^\sharp(|x|))^q |x|^{N(q-1)-q}dx           )^\frac{1}{q}\leq C\|D\chi_j(u_k)\|_{B_{kj}},
\end{eqnarray}
and
\begin{eqnarray}\label{2.11}
&&|B(0,1)|^\frac{q-1^\ast}{1^\ast q}(\int_{ A_{kj}       }((\chi_j(u_k))^\sharp(|x|))^q |x|^{N(q-1)-q}dx           )\nonumber\\
&\leq& C\|D(\chi_j(u_k))\|_{B_{kj}}(\int_{ A_{kj}       }((\chi_j(u_k))^\sharp(|x|))^q |x|^{N(q-1)-q}dx    )^\frac{q-1}{q},
\end{eqnarray}
where $B_{kj}=\{x\in \mathbb{R}^N \mid 2^{(n-1)(j-1)}\leq |u_k(x)|< 2^{(N-1)(j+2)} \}\supset A_{kj}$ and $\chi_{j_k}(u_k)=u_k$, $x\in A_{kj_k}$ by definition of $\chi_j(t)$. Let us sum up the above inequalities (\ref{2.11}) over $j\in \mathbb{Z}$. Note that by (\ref{1.17}). $\|D\chi_j(u_k)\|_{B_{kj}}\leq \|\chi'\|_{L^\infty}\|Du_k\|_{B_{kj}}$. Furthermore, one can break all the integers $j$ into four disjoint sets $J_1,J_2,J_3,J_4$,
\begin{eqnarray}\label{}
&&J_1=\{...,-4,0,4,8...\}\\
&&J_2=\{...,-3,1,5,9...\}\\
&&J_3=\{...,-2,2,6,10...\}\\
&&J_4=\{...,-1,3,7,11...\}
\end{eqnarray}
such that for any $m\in\{1,2,3,4\}$, all functions $\chi_j(u_k)$, $j\in J_m$, have pairwise disjoint supports. Consequently, $\sum_j\|Du_k\|_{B_{kj}}\leq4 \|Du_k\|$, we have therefore
\begin{eqnarray}\label{}
&&|B(0,1)|^\frac{q-1^\ast}{1^\ast q}(\int_{ \mathbb{R}^N       }(u^\sharp(|x|))^q |x|^{N(q-1)-q}dx           )\nonumber\\
&\leq& C\|Du_k\|\sup_{j}(\int_{ A_{kj}       }((\chi_j (u_k))^\sharp(|x|))^q |x|^{N(q-1)-q}dx    )^\frac{q-1}{q} \nonumber\\
&\leq& C\|Du_k\|\sup_{j}\|\chi_j (u_k)\|_{L^{1^\ast,q}(A_{kj})}^{q-1}.
\end{eqnarray}
It suffices now to show that for any sequence $\{j_k\}\subset \mathbb{Z}$, $\chi_{j_k}(u_k)\rightarrow 0$ in $L^{1^\ast,q}(A_{kj_k})$. Taking into account the invariance of the $L^{1^\ast,q}$-norm under the operators $g[j_k,y_k]$, it suffices to show that $\chi(2^{-j_k(N-1)}|u_k(2^{j_k}\cdot)|)\rightarrow 0$ in $L^{1^\ast,q}$, which is immediate by the assumption $g[j_k,y_k]u_k\rightharpoonup 0$, $\chi_{j_k}(u_k)=u_k$, $x\in A_{kj_k}$ and the argument of the \textbf{Step 1}, once we take into account for sequences $\chi(2^{-j_k(N-1)}|u_k(2^{j_k}\cdot)|)$ uniformly bounded in $L^\infty$. By $\chi\in C_0^\infty ((\frac{1}{2^{N-1}}, 4^{N-1}))$ and $\chi(t)=t$ as $t\in[1,2^{N-1}]$, we get $\chi(t)\leq t+1$ for all $t\in (\frac{1}{2^{N-1}}, 4^{N-1})$. Considering this fact and
\begin{eqnarray}\label{}
|u_k(x)|\leq 2^{(1+j_k)(N-1)},~x\in A_{kj_k},
\end{eqnarray}
we have
\begin{eqnarray}\label{}
\chi(2^{-j_k(N-1)}|u_k(2^{j_k}x)|)&\leq& 2^{-j_k(N-1)}|u_k(2^{j_k}x)|+1 \nonumber \\ \nonumber
&\leq& 2^{-j_k(N-1)}2^{(j_k+1)(N-1)}+1\\
&=&2^{(N-1)}+1,
\end{eqnarray}
which implies $\chi(2^{-j_k(N-1)}|u_k(2^{j_k}\cdot)|)$ bounded in $L^\infty$. The proof is completed.\\
\textbf{Theorem 2.3.} The imbedding $\dot{W}^{1,1}(\mathbb{R}^N)\hookrightarrow L^{1^\ast,q}(\mathbb{R}^N)$, $1<q\leq \infty$, $N\geq 2$, is cocompact with respect to the group $\mathcal{D}$ (\ref{1.16}).\\
\textbf{Theorem 2.4} (profile decomposition). Let $\{u_k\}\subset \dot{BV}(\mathbb{R}^N)$ be a bounded sequence. For each $n\in\mathbb{N}$ there exist $w^{(n)}\in \dot{BV}(\mathbb{R}^N)$ and sequences $\{j_k^n,y_k^n\}\subset \mathbb{Z}\times \mathbb{R}^N$ with $j_k^{(1)}=0$ and $y_k^{(1)}=0$ satisfying
\begin{eqnarray}\label{}
|j_k^{(n)}-j_k^{(m)}|+|y_k^{(n)}-y_k^{(m)}|\rightarrow\infty
\end{eqnarray}
whenever $m\neq n$, such that for a renumbered subsequence, $g[-j_k^{(n)},-y_k^{(n)}]u_k\rightharpoonup w^{(n)}$ as $k\rightarrow \infty$,
\begin{eqnarray}\label{}
r_k=u_k-\sum_{n}g[j_k^{(n)},y_k^{(n)}]w^{(n)}\rightarrow 0~\text{in}~L^{1^\ast,q}(\mathbb{R}^N),~q>1
\end{eqnarray}
where the series $\sum_n g[j_k^{(n)},y_k^{(n)}]w^{(n)}$ converges in $\dot{BV}(\mathbb{R}^N)$ uniformly in $k$, and
\begin{eqnarray}\label{}
\sum_n \|Dw^{(n)}\|+o(1)\leq \|Du_k\|\leq \sum_n \|Dw^{(n)}\|+\|Dr_k\|+0(1).
\end{eqnarray}
\textbf{Proof.} By Theorem 1.3 of \cite{AC} and Theorem 2.1, the claim easily follows.

\section{Counterexample}
In this section, we will show an example which indicates the imbedding $\dot{BV}(\mathbb{R}^N)\hookrightarrow L^{1^\ast,1}(\mathbb{R}^N)$ not cocompact. In fact we shall define a bounded sequence of $\dot{BV}(\mathbb{R}^N)$, which does not converge to zero in $L^{1^\ast,1}(\mathbb{R}^N)$. To this aim, fix a function
\begin{eqnarray}\label{}
\phi(x)=\left\{
  \begin{array}{ll}
    1,~&1< |x|\leq 2,~~~~~~~~~~~~~~~~~~~~~~~~~~~~~~~~~ \\
    0,~&|x|\leq 1~\text{or}~|x|>2,\\
  \end{array}
\right.
\end{eqnarray}
$x\in \mathbb{R}^N$, and obviously $\phi \in\dot{BV}(\mathbb{R}^N)$. We take the function $v_i$ given by the rescaling of $\phi$,
\begin{eqnarray}\label{}
v_i(x)=2^{i(N-1)}\phi(2^ix),~2^{-i}<|x|\leq2^{-(i-1)}.
\end{eqnarray}
For every positive integer number $n$, we choose $u_n$ as the following function:
\begin{eqnarray}\label{3.3}
u_n(x)=\frac{1}{n}\sum_{i=1}^n v_i(x).
\end{eqnarray}
Since the functions $v_i$ satisfies $\|v_i\|_{\dot{BV}}=\|\phi\|_{\dot{BV}}$ and (\ref{3.3}), we have
\begin{eqnarray}\label{}
\|u_n\|_{\dot{BV}}=\|\frac{1}{n}\sum_{i=1}^n v_i\|_{\dot{BV}}=\frac{1}{n}\sum_{i=1}^n\|v_i\|_{\dot{BV}}=\|\phi\|_{\dot{BV}}
\end{eqnarray}
and
\begin{eqnarray}\label{}
\|u_n\|_{L^{1^\ast}}^{1^\ast}=\frac{1}{n^{1^\ast-1}}\|\phi\|_{L^{1^\ast}}^{1^\ast}\rightarrow0,~ n\rightarrow\infty.
\end{eqnarray}
Then $u_n\rightarrow0$ strongly in $L^{1^\ast}(\mathbb{R}^N)$. It must be clearly $u_n\rightarrow0$ also in $L^{1^\ast,q}(\mathbb{R}^N)$, $1<q< 1^\ast$, as follows by interpolation since the sequence is bounded in $L^{1^\ast,1}(\mathbb{R}^N)$. We claim that $\{u_n\}$ does not converge to 0 in $L^{1^\ast,1}(\mathbb{R}^N)$ and therefore Theorem 2.1 does not hold with $q=1$. In fact, if $u_n\rightarrow 0$ strongly in $L^{1^\ast,1}(\mathbb{R}^N)$, then we have $u_n\rightharpoonup 0$ in $L^{1^\ast,1}(\mathbb{R}^N)=L^{\frac{N}{N-1},1}(\mathbb{R}^N)$. This implies for any linear functional $f\in (L^{1^\ast,1}(\mathbb{R}^N) )^\ast= (L^{\frac{N}{N-1},1}(\mathbb{R}^N))^\ast=L^{N,\infty}(\mathbb{R}^N)$, we have
\begin{eqnarray}\label{}
\lim_{n\rightarrow\infty}f(u_n)=0.
\end{eqnarray}
However, we can find a functional $f_0\in L^{N,\infty}(\mathbb{R}^N)$,
\begin{eqnarray}\label{3.7}
f_0(u)=\int_{\mathbb{R}^N}u(x)\frac{1}{|x|}dx.
\end{eqnarray}
Setting $u=u_n$ in (\ref{3.7}), we get
\begin{eqnarray}\label{3.8}
f_0(u_n)&=&\int_{\mathbb{R}^N}u_n(x)\frac{1}{|x|}dx\nonumber\\ \nonumber
&=&\frac{1}{n}\sum_{i=1}^n\int_{\mathbb{R}^N}v_i(x)\frac{1}{|x|}dx\\ \nonumber
&=&\frac{1}{n}\sum_{i=1}^n\int_{\mathbb{R}^N}2^{i(N-1)}\phi(2^ix)\frac{1}{|x|}dx\\ \nonumber
&=&\frac{1}{n}\sum_{i=1}^n2^{i(N-1)+i-Ni}\int_{\mathbb{R}^N}\phi(y)\frac{1}{|y|}dy\\ \nonumber
&=&\int_{\mathbb{R}^N}\phi(y)\frac{1}{|y|}dy\\ \nonumber
&=&\int_{\{y\in \mathbb{R}^N\mid 1<|y|\leq 2\}}\frac{1}{|y|}dy\\
&\geq&\frac{1}{2}|\{y\in \mathbb{R}^N\mid 1<|y|\leq 2\}|\neq0,
\end{eqnarray}
where $|A|$ is the $N$ dimensional Lebesgue measure of $A\subset\mathbb{R}^N$ in the last inequality. Taking $n\rightarrow\infty$ in the above inequality (\ref{3.8}), we have
\begin{eqnarray}\label{}
\lim_{n\rightarrow\infty}f_0(u_n)\geq\frac{1}{2}|\{y\in \mathbb{R}^N\mid 1<|y|\leq 2\}|\neq0.
\end{eqnarray}
This conclusion is a contradiction.

\end{document}